\documentclass[12pt]{amsart}

\usepackage[T2A,OT1]{fontenc}
\usepackage[ot2enc]{inputenc}
\usepackage[russian,english]{babel}

\usepackage{amssymb}
\usepackage{amscd}

\def\softl{l\kern-0.25ex\raise0.07ex\hbox{'}\kern-0.20ex}
\def\Drinfeld{{Drinfe\softl d }}

\newcommand{\injto}{\hookrightarrow}

\def\ZZ{{\mathbb{Z}}} 
 \def\Fp{{\mathbb{F}_p}} 
 \def\Fq{{\mathbb{F}_q}} 
\def\GG{{\mathbb{G}}} 
\def\PP{{\mathbb{P}}} 
\def\AA{{\mathbb{A}}} 

\def\Dc{{\mathcal{D}}} 
\def\Ec{{\mathcal{E}}} 
\def\Oc{{\mathcal{O}}} 

\def\Ell{{\mathcal{E}\ell\ell}} 

\def\C{{\mathbf{C}}}

\def\H{{\rm H}} 
\def\M{{\rm M}} 
\def\Spec{{\rm Spec}} 

\def\Lie{{\rm Lie}}

\def\rk{{\rm rk}}

\def\op{{\rm op}}

\def\deg{{\rm deg}}

\def\an{{\rm an}}

\def\pr{{\rm pr}}
\def\i{{\imath}}

\swapnumbers
\theoremstyle{plain}
\newtheorem{theorem}[subsection]{Theorem}
\newtheorem{corollary}[subsection]{Corollary}

\newtheorem{proposition}[subsection]{Proposition}
\newtheorem{definition}[subsection]{Definition}
\theoremstyle{definition}
\newtheorem{example}[subsection]{Example}

\newtheorem{void}[subsection]{}

\begin{document}

\markboth{Lenny Taelman}
 {$\Dc$-elliptic sheaves and uniformisation}

\title
{$\Dc$-elliptic sheaves and uniformisation}

\author{Lenny Taelman}
\email{lenny@math.rug.nl}

\maketitle

%
%
%

\thispagestyle{empty}

\section{Introduction}

This text seeks to expose a function field analogue - in the
sense of $\Fp[t]$ versus $\ZZ$ - to 
objects such as modular curves and quaternionic Shimura curves.

Over the field of complex numbers, these curves permit a dual
description as either algebraic varieties or Riemann surfaces. The
former is obtained as a parameter space for certain Abelian varieties,
the latter as a quotient space of a homogeneous space, to wit the complex upper half plane, by a discrete group action. 

Modelled after the classical modular curves are the \Drinfeld modular
curves and the higher dimensional \Drinfeld modular varieties.
Similarly, they can be treated either algebraically or analytically  (see \cite{Drinfeld74}). A wider class of
varieties, containing at once the \Drinfeld modular varieties and the
proper counterparts of the quaternionic Shimura curves, consists of
the moduli spaces of Laumon, Rapoport and Stuhler. They are
algebraic varieties that parametrise certain algebraic objects,
namely $\Dc$-elliptic sheaves (see \cite{Laumon93}). However,
the dual, analytic, description is absent from the literature.

Thus, in a slightly more precise formulation, this manuscript
aims to give the moduli spaces for $\Dc$-elliptic sheaves an
analytic description as quotients of 
homogeneous spaces. This proceeds roughly in three steps. First,
$\Dc$-elliptic sheaves are shown to form a particular subclass of the
category of Anderson $t$-motives (\S 2, proposition
\ref{equivalence}). Second, these Anderson $t$-motives
are proven to be \emph{uniformisable}, meaning that they can be
obtained as quotients of vector spaces by lattices, and the class of
lattices so obtained is determined (\S 3,
propositions \ref{uniformizability} and \ref{fullness}). Finally, these
lattices are parametrised by quotients of homogeneous spaces; this
yields the main result (\S 4, theorem \ref{thm}).

This uniformisation of the varieties of Laumon, Rapoport and Stuhler
at the ``infinite'' valuation complements recent work of Hausberger
(\cite{Hausberger05}). In the spirit of \Drinfeld and Cherednik's
$p$-adic uniformisation of Shimura curves, he treats the
analytic structure at the ``ramified'' places. Hausberger has confided
to the author that he has also obtained, but not published, a
uniformisation at the ``infinite'' valuation.

\section{$\Dc$-elliptic sheaves and Anderson $A$-motives}

\begin{void}
\marginpar{\tiny standard Drinfeld data} 
The following data are fixed. A complete, smooth and
geometrically irreducible curve $X$ over a finite field $\Fq$ of
characteristic $p$ and a closed
point $\infty\in X$ called infinity. The degree of the residue
class field of $\infty$ over $\Fq$ is denoted by $\alpha$. The
function field of $X$ is denoted by $F=\Fq(X)$ and the subring of
functions which are regular
everywhere except possibly at infinity by $A=\H^0(X-\infty,\Oc_X)$.
On $A$, $\deg$ is defined as the pole order at infinity.
The completion of $F$ with respect to the infinite valuation is
$F_\infty$. Also, a complete and algebraically closed field $\C$
containing $F_\infty$ is fixed and the resulting embeddings denoted by
$\i:F_\infty\to \C$, $\i:F\to\C$ and $\i:A\to\C$. Of course, the kernel
of $\i:A\to\C$ is trivial, this corresponds to what is usually called
``infinite characteristic'', and reflects the particular interest into
analytic structures taken in this exposition. Put differently:
reductions modulo primes of $A$ will not be touched upon. Finally,
denote by $\hat{A}$ the product of all completions $A_v$ over all
places $v\neq \infty$ and by
$\AA_A=\hat{A}\otimes_A F$ the ad\`{e}le ring away from
infinity.
\end{void}

\begin{void}
\marginpar{\tiny the central simple algebra} 
Also, a central simple $F$-algebra $D$, unramified at $\infty$ is fixed.
Denote the dimension of $D$ over its centre by $d^2$. Let $L$ be an 
algebraically closed field containing $\Fq$, then 
$D\otimes_{\Fq}L\cong D\otimes_F L(X)$ is isomorphic to the full matrix algebra $\M(d,L(X))$, since $L(X)$ is a $C_1$ field.

Fix a subsheaf of $\Oc_X$-algebras $\Dc$ inside the constant sheaf
$D$ with the property that for all points $x\in X$, $\Dc_x$ is
a maximal order in $D_x$ (it is shown in \cite{Laumon93} that such sheaves exist.) Then $O_D:=\H^0(X-\infty,\Dc)$ is a maximal
$A$-order in $D$. 
\end{void}

\begin{void}
\marginpar{\tiny the sheaves} 
Let $B$ be a $\C$-algebra and denote $X\times_{\Fq}\Spec(B)$ by
$X_B$. The curve $X_B$ on $\Spec(B)$ comes equipped with a
section $\i^\#:\Spec(B)\to X_B$ which derives from the
embedding $\i$ of the function field of $X$ into $\C$ and takes as
constant value the generic point of $X$. 

The projections $X_B \to X$ and $X_B\to \Spec(B)$ are denoted
by $\pr_X$ and $\pr_S$ respectively. The
pullback of the geometric $q$-Frobenius on $\Spec(B)$ to
$X_B$ by $\pr_X$ by $\sigma$. Thus $\sigma$ acts trivially on
$X$ and as the $q$-th power endomorphism on
the ring of functions $B$.
\end{void}

\begin{definition} \label{ellsh}
 A \emph{$\Dc$-elliptic sheaf} on $X_B$ is a commutative
diagram:
\[ \begin{CD}
\ldots @>j>> \Ec_{-1} @>j>> \Ec_0 @>j>> \Ec_1 @>j>> \Ec_2 @>j>> \ldots \\
  @.    @AtAA @AtAA @AtAA @AtAA \\
\ldots @>j>> \sigma^\ast\Ec_{-2} @>j>> \sigma^\ast\Ec_{-1} @>j>>
\sigma^\ast\Ec_0 @>j>> \sigma^\ast\Ec_1 @>j>>
\ldots \end{CD} \]
where the $\Ec_i$ are locally free right $\pr_X^\ast\Dc$-modules of rank
$1$ and  the $t$ and the $j$ are $\pr_X^\ast\Dc$-linear
injections. These data are constraint to the following conditions:
\begin{enumerate}
 \item \emph{Periodicity}: the composition of $d\alpha$
  consecutive $j$'s defines an identification
  $\Ec_{i+d\alpha}=\Ec_i\otimes_{\Oc_{X_B}} \
  (\Oc_X(\infty)\boxtimes \Oc_{\Spec(B)})$
  through the natural embedding $\Oc_X\injto \Oc_X(\infty)$;
 \item \emph{Pole}: the direct image of $\Ec_i/j(\Ec_{i-1})$ by
  $\pr_B$ is a locally free $\Oc_{\Spec(B)}$-module of rank
   $d$;
 \item \emph{Zero}: $\Ec_i/t(\sigma^\ast\Ec_{i-1})$ is 
 the direct image by $\i^\#:\Spec(B)\to X_B$ of a locally free
 $\Oc_{\Spec(B)}$-module.
\end{enumerate}
 A shorthand notation for such an object is $(\Ec_i,j,t)$. A
 \emph{morphism of $\Dc$-elliptic sheaves} on $X_B$ is a
 collection of $\pr_X^\ast\Dc$-linear morphisms $\Ec_i\to\Ec'_i$ 
 that commute with the $j$'s and the $t$'s.
\end{definition}

\begin{void}\label{decalage}
\marginpar{\tiny d\'{e}calage} 
  It will be necessary to consider elliptic
 $\Dc$-sheaves up to \emph{d\'{e}calage}. Formally this can be done by
 allowing more morphisms. Namely, a morphism of  $\Dc$-elliptic sheaves up to
 d\'{e}calage consists of a fixed integer $n$ and a collection
 of $\pr_X^\ast\Dc$-linear morphisms
 $\Ec_i\to\Ec'_{i+n}$ commuting with the $j$'s and $t$'s.
 \end{void}

\begin{definition} 
\marginpar{\tiny the modules} 
 An Anderson $A$-motive with multiplication by $O_D$ over $\C$
 is a left $O_D^\op\otimes_\Fq \C[\sigma]$-module $M$ satisfying
 \begin{enumerate}
  \item $M$ is projective of rank $1$ over $O_D^\op\otimes \C$
  \item $M$ is projective of rank $d$ over $\C[\sigma]$
  \item there exists a positive integer $n$ such that
    $(a-\i(a))^n(M/\sigma M)=0$ for all $a\in A$.
 \end{enumerate}
\end{definition}

\begin{void}\label{altdef}
 Denote by $\sigma$ the endomorphism of $O_D^\op\otimes_\Fq \C$
 that acts trivially on $O_D^\op$ and as the $q$-th power map
 on $\C$.
 Then an Anderson $A$-motive with multiplication by $O_D$ can
 also be described as a
 $O_D^\op\otimes \C$-module $M$ on which the action of $\sigma$ is
 described by either a linear morphism
 $\tau:\sigma^\ast M \to M$ or a semi-linear map from $M$ to
 itself, also (abusively) denoted by $\tau:M\to M$.
\end{void}

 As is already suggested in \cite[\S 3]{Laumon93}, there is a close
 relationship between $\Dc$-elliptic sheaves and Anderson
 $t$-motives (in this manuscript called $A$-motives). This
 kinship is best understood through Drinfe\softl d's
 theory of vector bundles on the non-commutative projective line.
 For more details on these objects, one can consult \S 3
 of \emph{loc. cit.}, which establishes most of what is needed
 in order to properly relate Anderson $A$-motives with
 multiplication by $O_D$ to $\Dc$-elliptic sheaves. In the following
 paragraphs the remaining steps are taken.

\begin{void}\label{begineq}
\marginpar{\tiny from sheaves to modules} 
 Since by the pole axiom, the subsequent quotients
 $\Ec_{i+1}/\Ec_i$ in
 a $\Dc$-elliptic sheaf $(\Ec_i,j,t)$ on $X_\C$ are
 supported at $\infty$,
 the $O_D^\op\otimes_k \C$-module $M:=\H^0(X_\C-\infty,\Ec_i)$ is
 independent of $i$. Therefore the
 $t:\sigma^\ast \Ec_i \to \Ec_{i+1}$ induce a morphism 
 $\tau:\sigma^\ast M\to M$ which turns the
 $O_D^\op\otimes_k \C$-module
 $M$ into a $O_D^\op\otimes_k \C[\sigma]$-module
 (as in \ref{altdef}). By construction,
 $M$ depends only on the given $\Dc$-elliptic sheaf up to
 d\'{e}calage. It follows from \cite[3.17]{Laumon93} that $M$
 is an Anderson $A$-motive with multiplication by $O_D$.
\end{void}

\begin{void}
\marginpar{\tiny from modules to sheaves} 
 The above construction is functorial and it is possible to
 define an inverse functor. It suffices to extend the given
 $\C[\sigma]$-module $M$ to a vector bundle on the non-commutative
 projective line with coordinate $\sigma$ and then apply the
 equivalence 3.17 of \emph{loc. cit.} This is tantamount to
 a module-theoretic construction, that will be described next,
 first in the rank one \Drinfeld module case, that is, for $D=F$,
 and then for a general $D$.
\end{void}

\begin{void}
 Let $D=F$. The non-commutative ring of power series in
 $\sigma^{-1}$, denoted by $\C[[\sigma^{-1}]]$, has the skew field
 of Laurent series $\C((\sigma^{-1}))$ as quotient field. Given
 an $M$, a Anderson $A$-motive with multiplication by
 $O_D=A$, under the aforementioned equivalence, the definition
 of a $\Dc$-elliptic sheaf $(\Ec_i,j,t)$ such
 that $M=\H^0(X_\C-\infty,\Ec_i)$ boils down to the definition of a
 free rank one $\C[[\sigma^{-1}]]$-submodule
  $$ W \subset \C((\sigma^{-1})) \otimes_{\C[\sigma]} M .$$
 The right-hand side of this inclusion is isomorphic to
 $\C((\sigma^{-1}))$, and the corresponding submodules are precisely
 those generated by a $\sigma^m$ for an integer $m$. It can be
 checked that every such submodule defines a vector bundle on
 the non-commutative projective line with the correct
 properties, hence also an elliptic sheaf. Different
 choices of $m$ define elliptic sheaves that are
 equivalent under d\'{e}calage.
\end{void}

\begin{void}\label{morita}
\marginpar{\tiny Morita equivalence} 
 Extending the above reasoning to the general case is a matter
 of applying Morita equivalence. As this equivalence
 is also central to several of the proofs yet to come, it is
 worthwhile to recall a precise formulation here.
 Let $R$ be a ring with unity. Then to every $R$-module $N$ it is
 possible to associate the $\M(n,R)$-module $N\oplus N\oplus \cdots
 \oplus N$, the direct sum of $n$ copies of $N$.
 \emph{Morita equivalence} states that this
 defines an equivalence from the category of  $R$-modules to
 the category of $\M(n,R)$-modules, and that this equivalences
 preserves such properties as projectivity, injectivity, and
 finite generation. Also, the construction is functorial in $R$.
 A detailed exposition can be found in \cite[\S 7]{Lam99}.
\end{void}

\begin{void}\label{endeq}
\marginpar{\tiny from modules to sheaves (continued)}
 The restriction on $D$ is now dismissed. Now
 $\C((\sigma^{-1})) \otimes_{\C[\sigma]} M$ is a module over
 $F_\infty \otimes_A O_D^\op$ and the submodule
  $$ W \subset \C((\sigma^{-1})) \otimes_{\C[\sigma]} M $$
 is required to be isomorphic to $\C[[\sigma^{-1}]]^d$ and to
 be equipped with a continuous action by
 $\Dc_\infty\subset F_\infty\otimes_A O_D^\op$. Since $\Dc_\infty$
 is a matrix algebra inside the matrix algebra
 $F_\infty\otimes_A O_D^\op$, the machinery of Morita applies.
 It turns the set of such $W$ into the set of
 $\C[[\sigma^{-1}]]$-submodules of $\C((\sigma^{-1}))$. As above,
 these are generated by some $\sigma^m$ and different $m$ define
 different representatives of the same class of 
 $\Dc$-elliptic sheaves up to d\'{e}calage.
\end{void}

The above discussion (\ref{begineq}-\ref{endeq}) sums up to:

\begin{proposition}\label{equivalence}
 The category of $\Dc$-elliptic sheaves on $X_\C$ up to
 d\'{e}calage is equivalent to the category of
 Anderson $A$-motives with multiplication by $O_D$.
\end{proposition}

\begin{example}
\marginpar{\tiny \Drinfeld modules} 
 If $O_D$ is the algebra of $d\times d$ matrices over $A$ then Morita
equivalence implies the equivalence of Anderson $A$-motives with
multiplication by $O_D$ and Anderson $A$-motives of dimension $1$ and
rank $d$. Thus the theory of $\Dc$-elliptic sheaves includes the
\Drinfeld modules as a special case.
\end{example}

\begin{void} 
\marginpar{\tiny level structures} 
Given an effective divisor $I$ on $X-\infty$ it is possible to define
level $I$-structures on $A$-motives with multiplication by $O_D$
as well
as on $\Dc$-elliptic sheaves such that the equivalence between the two
categories extends to an equivalence of the categories ``with level
structure''. For the sake of clarity, the results in this text
are first stated and proved ignoring level structures, and
the corresponding propositions with level structure are simply
stated as corollaries. The omitted proofs are never hard.
\end{void}

The $I$ in the succeeding definitions refers in \emph{abus de
notation} to both the divisor $I$ on $X$ and the corresponding
ideal $I\subset A$. 

\begin{void} 
\marginpar{\tiny level structures on sheaves} 
Since $\infty$ and $I$  are disjoint, the
restriction $\Ec_I$ of the $X_B$-sheaf $\Ec_i$ to $I_B$ is
independent of $i$ and for the same reason the injection $t$ restricts
to an isomorphism $t:\sigma^\ast\Ec_I\to \Ec_I$.  
A \emph{level $I$-structure on a $\Dc$-elliptic sheaf $(\Ec_i,j,t)$}
is an isomorphism of $\Oc_{I_B}$-sheaves:
$\Dc_I\boxtimes\Oc_B \to \Ec_I$ such
that the Frobenius $\sigma^\ast \Oc_B\to \Oc_B$ corresponds to 
$t:\sigma^\ast\Ec_I\to \Ec_I$.
\end{void}

\begin{void} 
\marginpar{\tiny level structures on modules} 
 A \emph{level $I$ structure on an Anderson $A$-motive
 with multiplication by $O_D$} is the choice of an isomorphism
 $(O_D^\op/IO_D^\op)\otimes_{\Fq} \C\to M/IM$
 such that the Frobenius on $\C$ matches the action of $\sigma$
 on $M/IM$.
\end{void}

\section{Uniformisation of the Anderson $A$-motives}

\marginpar{\tiny introduction} 
Whereas all Abelian varieties over the field of the complex numbers
can be obtained as the quotient of a complex vector space by a
lattice, the parallel statement for Anderson $A$-motives is false
in general. Those Anderson $A$-motives that \emph{can} be obtained
in such a way are called  \emph{uniformisable}. It 
will be shown in this section that the $A$-motives considered above,
\emph{videlicet} the $A$-motives with multiplication by $O_D$, are
in fact uniformisable.

\begin{void}
\marginpar{\tiny on uniformisability} 
 A very brief overview of the results of \cite{Anderson86}
 is given. There is an anti-equivalence of categories that
 assigns to an Anderson $A$-motive $M$ of dimension $d$,
 an action of $A$ on an algebraic group $E(M)\cong\GG_a^d$,
 satisfying certain properties. This $E(M)$ occurs naturally
 in an exact sequence of \emph{analytic} $A$-modules:
\[ \begin{CD}
  \Lambda @>>> \Lie(E(M)) @>{\exp_{E(M)}}>> E(M)
  \end{CD} \]
 Here $\Lie(E(M))$ is a $d$-dimensional vector space over $\C$
 on which $A$ acts by the embedding $\i:A\injto \C$, and $\Lambda$
 is a lattice, \emph{i.e.} a discrete and projective $A$-submodule.
 The Anderson $A$-motive $M$ can be recovered from this lattice
 if $\exp$ is surjective. There exist, however, $M$ for which
 this is not the case. $M$ is called \emph{uniformisable}
 when $\exp_{E(M)}$ \emph{is} surjective and by
 \cite[\S 2]{Anderson86} this
 is the case precisely when the rank of $\Lambda$ equals the rank
 of $M$ (in general $\rk_A(\Lambda)\leq \rk_{A\otimes\C}(M)$). \end{void}

\begin{void}
\marginpar{\tiny a criterion for uniformisability}
 Fix a subring $\Fp[a]\subset A$ such that $a$ has a pole of order
 prime to $p$ at $\infty$ and such that the resulting fraction
 field extension $F/\Fp(a)$ is separable. If the rank of $A$
 over $\Fp[a]$ is
 $r$, then an Anderson $A$-motive with multiplication by $O_D$ is an
 Anderson $\Fp[a]$-motive of rank $rd^2$ and dimension $d$. Given an
 Anderson $A$-motive over $\C$ it is now possible to extend its
 scalars from $\C[a]$ to $\C\{\{a\}\}$ - the ring of power series with
 convergence radius at least $1$ - which yields a module
 $\tilde{M}=M\otimes_{\C[a]}\C\{\{a\}\}$. The semilinear action
 $\tau:M\to M$ extends to a semilinear
 $\tilde{\tau}:\tilde{M}\to\tilde{M}$ and $M$ is said to be
 \emph{analytically trivial} if $\tilde{M}$ has a $\C\{\{a\}\}$-basis
 consisting of $\tilde{\tau}$-invariant vectors. This is all
 that is needed to
 state Anderson's criterion for uniformisability
 \cite[\S 2]{Anderson86}. Namely, it states that $M$ is
 uniformisable if and only if it is analytically trivial. 
\end{void}

\begin{proposition}\label{uniformizability}
\marginpar{\tiny uniformisability of the modules}
 Over $\C$, all $d$-dimensional $A$-motives with multiplication
 by $O_D$ are uniformisable.
\end{proposition}

\begin{proof}
 Fix an isomorphism
 $O_D^\op\otimes_{\Fq}\C\cong \M(d,A\otimes_\Fq\C)$.

 Let $M$ be  an $A$-motive as in the proposition. Then $M$ is a
 projective $O_D^\op\otimes \C$-module of rank one together with a 
 semilinear map $\tau:M\to M$ commuting with the action
 of $O_D^\op$
 and having a determinant equal to $(a-\i(a))^d$ up to a unit.
 Morita equivalence associates to these data a module $N$ of rank
 $n/d$ over $A\otimes \C$ equipped
 with a semilinear $\tau_N:  N\to N$ of determinant a
 unit times $(a-\i(a))$ such that $M$ can be recovered from $N$
 as a $d$-fold direct sum: $M\cong N\oplus N \oplus \cdots \oplus N$. 
 The $A$-motive $N$ has dimension $1$, thus determines a
 \Drinfeld module. \Drinfeld modules are uniformisable
 (\cite{Drinfeld74}), and it
 follows from Anderson's analytic triviality criterion that $N$ is
 analytically trivial. This shows that
 $\tilde{M}=\tilde{N}\oplus\tilde{N}\oplus\cdots\oplus\tilde{N}$
 has a $\tilde{\tau}$-invariant $\C\{\{a\}\}$-basis, hence $M$
 is uniformisable.
\end{proof}

\begin{void}
\marginpar{\tiny lattices and representations} 
Thus, over $\C$, to every $d$-dimensional $A$-motive $M$ with
multiplication by $O_D$, or to every $\Dc$-elliptic sheaf on
$X_\C$, a discrete lattice $\Lambda$ of $A$-rank $r^2$ in
$V:=\Lie(E(M))$ is assigned. This assignment is a faithful functor
on the category of uniformisable $A$-motives
\cite[2.12.2]{Anderson86}, and this immediately implies that
$\Lambda\subset V$ is more than just an $A$-lattice in a
$\C$-vector space: it naturally carries the structure of an $O_D$-lattice
in a faithful $\C$-linear representation of $D$. The following proposition 
confirms that every pair $\Lambda\subset V$ of a rank one
$O_D$-lattice in a faithful
$d$-dimensional representation of $D$ occurs in such a way.
\end{void}

\begin{proposition}\label{fullness}
\marginpar{\tiny surjectivity of the construction} 
 Let $\Lambda\subset V$ be a discrete projective
 and rank one $O_D$-submodule of a $d$-dimensional faithful linear
 representation of $D$. There exists a uniformisable Anderson $
 A$-motive with multiplication by $O_D$ and
 an isomorphism $\Lie(E(M))\to V$ that surjects $\ker(\exp_{E(M)})$
 onto $\Lambda$.
\end{proposition}

\begin{proof}
 By \cite[corollary 3.5.1]{Anderson86}, it is sufficient to find
 an $A$-lattice
 $\Lambda'\subset V$ whose $F_\infty$-span equals the $F_\infty$-span of
$\Lambda$ and which is the kernel of the exponential function of some
uniformisable abelian $A$-module. Since $D$ is unramified at the infinite place,
$D\otimes F_\infty$ is isomorphic to $\M(d,F_\infty)$, hence the
$F_\infty$-span of $\Lambda$ splits under the action of the diagonal
idempotents as a direct sum of $d$ terms each lying within a one-dimensional linear
subspace of $V$. Each such term is the span of a suitable
lattice in a one-dimensional linear space and the direct sum of these
lattices is the $\Lambda'$ which is being sought for. It is the kernel
of the exponential function of an abelian $A$-module which is
uniformisable because it is the
direct sum of $d$ suitable \Drinfeld modules.
\end{proof}

\begin{corollary}\label{eq2}
 The category of Anderson $A$-motives with multiplication by $O_D$
 is anti-equivalent with the category of pairs $(V,\Lambda)$, where
 $V$ is a faithful $d$-dimensional representation of $D$ and
 $\Lambda\subset V$ an $O_D$-lattice.
\end{corollary}

\section{Analytic moduli spaces of $\Dc$-elliptic sheaves}

In this section, all $\Dc$-elliptic sheaves are considered up to
d\'{e}calage (see \ref{decalage}).

The requisite preparations have now been made in order to
give the varieties of Laumon, Rapoport and Stuhler an analytic
description. The main ingredient is the classification of the 
pairs $(V,\Lambda)$ of corollary \ref{eq2}. These are classified
by quotients of the \Drinfeld symmetric space
$\Omega^d(\C)=\PP^{d-1}(\C)-\PP^{d-1}(F_\infty)$, in very much
the same way that some Shimura varieties are constructed as
quotients of the complex upper half plane. In what follows
$D^\ast(R)$ denotes $(D\otimes_F R)^\ast$, for an $F$-algebra $R$.

\begin{proposition}\label{classlat}
\marginpar{\tiny classification of lattices}
 There is a natural bijection between the set of isomorphism
 classes of pairs $(V,\Lambda)$ of $d$-dimensional $\C$-linear
 representations of $D$ with $O_D$-lattice $\Lambda\subset V$ and
 the double co-sets in
  \[ D^\ast \backslash
   \left[\Omega^{d}(\C) \times D^\ast(\AA_A)  /
    U \right] \]
  where $D^\ast$ acts diagonally on the product and in particular on
  $\Omega^{d}(\C)$ by the choice of an isomorphism
  $D\otimes_F F_\infty\cong \M(d,F_\infty)$ and where $U$ is
  the compact open subgroup
  $(O_D\otimes_A\hat{A})^\ast\subset D^\ast(\AA_A)$.

\end{proposition}

\begin{proof}
 Start of with a pair $(V,\Lambda)$. Since all
 representations are conjugate, it is possible to assume without
 loss of generality that
 $V=F_\infty^d\otimes \C=\C^d$ on which $D$ acts by a fixed
 isomorphism $D\otimes_A F_\infty \cong \M(d,F_\infty)$.
 Consider the $F$-span
 $F\Lambda$ of the lattice. This is a free module of rank $1$ over
 $D$ lying inside $V$ and the choice of a generator
 marks a point on $\PP^{d-1}(\C)$ and identifies $F\Lambda$ with
 $D$. The marked point lies in
 $\Omega^d(\C)\subset\PP^{d-1}(\C)$ by the discreteness of $\Lambda$.
 The embedding $\Lambda\subset D$ can be tensored to an
 embedding $\hat{A}\Lambda \subset D\otimes \AA_A$
 and the former can be recovered from the latter as
 $\Lambda=\hat{A}\Lambda \cap D\otimes 1$. But all projective
 $O_D\otimes\hat{A}$-modules are free, consequently the projective
 $O_D$-submodules $\Lambda\subset D$ of rank one are in
 bijection with the free rank
 one $O_D\otimes\hat{A}$-submodules of $D\otimes\AA_A$ and the latter
 are in bijection with
 $(D\otimes\AA_A)^\ast/(O_D\otimes\hat{A})^\ast$. It remains to
 mod out by the choice of the generator of $F\Lambda$, \emph{i.e.}
 by $D^\ast$, to establish the desired one-to-one
 correspondence. 
\end{proof}

\begin{void}
\marginpar{\tiny analysis of the proposition}
Chaining the last proposition with the equivalences \ref{eq2} and
\ref{equivalence} furnishes a bijection between the set
of isomorphism classes of $\Dc$-elliptic sheaves on $X_\C$
and the above double co-sets. This double co-set space inherits the
structure of a rigid analytic space from the rigid analytic
space $\Omega^d$ (see \cite{Drinfeld74}). So, at least point-wise,
$\Dc$-elliptic sheaves are classified by a rigid analytic space, but
this set-theoretic result is not enough to give an analytic structure
to the algebraic moduli spaces of Laumon, Rapoport and Stuhler.
\end{void}

\begin{void}
\marginpar{\tiny analytic moduli} 
A close examination of the proof of proposition \ref{classlat}
shows that it establishes a slightly stronger result. In fact,
it demonstrates that every analytic family of pairs $(V,\Lambda)$
over a rigid analytic $\C$-space $Y$ - the definition of such a 
family should be clear - results in an analytic map from
$Y$ to the double co-set space. In fact, also this stronger statement
can be easily pulled through the chain of equivalences \ref{eq2} and
\ref{equivalence}, using the results of \cite{Boeckle05}. The
outcome is a procedure that assigns to every $\Dc$-elliptic sheaf
on $X_B$ for some $\C$-algebra $B$, a morphism of rigid
analytic $\C$-spaces from $\Spec(B)^\an$ to the double co-set space
of the proposition. 
\end{void}

\begin{void}
\marginpar{\tiny algebraic moduli} 
 In \cite{Laumon93} it is shown that when $I\neq A$, 
 the functor that ``maps'' a $\C$-algebra $B$ to the set of
 isomorphism classes of $\Dc$-elliptic sheaves on $X_B$ is
 representable by a quasi-projective $\C$-scheme of
 dimension $d-1$, denoted by $\Ell_{X,\Dc,I}$. Taking the
 quotient by a finite group yields a coarse moduli
 scheme $\Ell_{X,\Dc,A}$ classifying 
 $\Dc$-elliptic sheaves without level structure. All these
 varieties are shown to be smooth and under the condition that
 $D$ be a division algebra they are shown to be complete
 in \emph{loc.cit.}
\end{void}

\begin{void}
\marginpar{\tiny relating algebraic and analytic moduli}
  Now the above considerations define a natural map of rigid
  analytic spaces
   \[ \Ell_{X,\Dc,A}(\C)^\an \to D^\ast \backslash
   \left[\Omega^{d}(\C) \times D^\ast(\AA_A)  /
    U\right], \]
  which is a bijection on the sets of $\C$-valued points. Since both
  spaces are reduced, this has to be an isomorphism. This proves the
  fundamental case of the main theorem.
\end{void}

\begin{theorem}
\marginpar{\tiny principal case of the theorem}
 There is a natural isomorphism of rigid analytic spaces
   \[ \Ell_{X,\Dc,A}(\C)^\an \cong D^\ast \backslash
   \left[\Omega^{d}(\C) \times D^\ast(\AA_A)  /
    U\right]. \]
\end{theorem}

\begin{void}
\marginpar{\tiny free vs. projective} 
 A triple $(\Ec_i,j,t)$ is said to be \emph{free} if and only
 if the restrictions of the $\Ec_i$ to $X\times s$ are free for
 all geometric points $s$ of the base. For example, when $X$ is the
 projective line and $\infty$ is rational then all elliptic
 $\Dc$-sheaves are free. This notion
 readily translates to Anderson $A$-motives. $M$ being free then
 signifies $M$ being free as $A\otimes \C$-module.
 Free $\Dc$-elliptic sheaves are
 classified by closed and open subspaces
 $\Ell^0_{X,\Dc,I}\subset \Ell_{X,\Dc,I}$. 
  The above theorem
 quite easily generalises to the main result of this exposition.
\end{void}

\begin{theorem}\label{thm}
\marginpar{\tiny main theorem} 
 The analytification of the moduli spaces $\Ell^0$ and $\Ell$ is
 $$ \Ell^0_{X,\Dc,I}(\C)^\an \cong G(I) \backslash \Omega^{d}(\C) $$
and
 $$ \Ell_{X,\Dc,I}(\C)^\an \cong
   D^\ast \backslash \left[\Omega^{d}(\C)
   \times D^\ast(\AA_A) / U(I)\right]
$$
where $G(I)\subset O_D^\ast$ and $U(I)\subset U$ are the subgroups
of elements
that reduce to $1$ modulo $I$ and where $D^\ast$ acts on
$\Omega^{d}(\C)$ by the choice of an isomorphism
$D\otimes_F F_\infty\cong \M(d,F_\infty)$.
\end{theorem}

If $D=\M(d,F)$, then the above result boils down to the
uniformisation of the \Drinfeld modular varieties $M^d_I$
as established in \cite{Drinfeld74}. 

\subsection{Example}
\marginpar{\tiny quaternionic Shimura curves}
When $O_D$ is a maximal order in a quaternion algebra over $F$ that is
unramified at $\infty$, then for every level $I$, the moduli space described
above is a complete smooth curve over $\C$. By the above theorem, these spaces
are quotients of $\Omega^2\subset \PP^1(\C)$ by some subgroup $G$ of
$O_D^\ast\subset D\otimes_F \C\cong \M(2,\C)$. It follows from
the theory of Mumford curves (see \emph{e.g.}
\cite{Mumford72}, \cite{Gerritzen80}) that the genus of this curve is
the rank of $G$, \emph{i.e.} the rank of the Abelianisation of $G$. It
is in general however not straightforward to calculate the rank of a
given $G$. Usually it is easier to find a larger group $H \supset G$
such that $H\backslash\Omega^2$ has genus zero, and to calculate the
genus of $G\backslash\Omega^2$ from the ramification of
$G\backslash\Omega^2\to H\backslash\Omega^2\cong\PP^1$.

\bibliographystyle{plain}
\bibliography{master}

\end{document}